# Parameterizing Fullerenes with Vertex Combinations

Li, Shaoqing

Anhui Jianzhu University, Hefei, Anhui, China.

**Abstract:** A fullerene graph can be embedded in a piecewise linear 2-manifold with each non-hexagonal carbon ring corresponding to a cone vertex. Adjacent two or three such vertices can be combined as a cluster cut out from a parent cone round a parent vertex. The locations for the combined vertices are distributed on their parent cone with certain symmetries. With proper combinations, a fullerene can be simplified to a parent structure easily to be parameterized. And then, with additional parameters for the combined vertices in each cluster, the fullerene can also be parameterized. Such parameters may be different for fullerenes from type to type, but they are able to characterize the specific shape for each type.

## 1. Introduction

A fullerene graph (fullerene for short) is a polyhedral graph modeling an allotrope of carbon tiling with pentagonal and hexagonal carbon rings in symmetries. Researchers have focused on the fullerenes for their symmetries and parameterization ever since the discovery of the allotrope in 1985 [1-3]. With icosahedrons (not necessarily regular) as their master polyhedrons, the present author has studied their intrinsic symmetries and accomplished a general parameterization. However, in some cases, such as carbon nanotubes, the icosahedron is not a perfect model as the parameters cannot directly characterize their shape and their values may be negative [4].

A fullerene can be embedded in a piecewise linear 2-manifold with each non-hexagonal carbon ring corresponding to a cone vertex [4, 5]. Distinct from the vertices of fullerene graph, the vertices mentioned in this paper are all cone vertices of the fullerene manifold, and the edge number of the carbon ring is called its degree.

On the fullerene manifold, not only the area around an vertex is a cone manifold, but the area around a cluster of vertices is also part of a master cone [1], or cylinder in special [4]. Adjacent vertices can be combined as a cluster cut out from the master cone. If so, the master cone is called the parent cone and its virtual vertex is called the parent vertex of the cluster.

The parent vertex is determined by the combined vertices. Therefore, around the parent cone, the combined vertices are not independent. There are certain symmetries of combined vertices around their parent cone. With these symmetries, the combined vertices can be parameterized.

With appropriate vertex combinations, a fullerene may have a simple parent structure, such as (3,6)-fullerene or (4,6)-fullerene, that can be easily parameterized [4]. If so, with the additional parameters for the combined vertices, the fullerene can be also parameterized.

In 1998, William P. Thurston mentioned that fullerenes (just as their duals) can be parameterized with manifold developing, i.e. vertex combination [1]. However, the study has not attracted enough attention and the related exploration has been delayed for decades. In this paper, vertex combinations were studied and parameterizations of fullerenes with different shapes were explored in light of the combinations.

The remainder of the paper is as follows. In section 2 and section 3, the combination of two vertices and the combination of three vertices were studied respectively, and their symmetries around the parent cone were revealed. In section 4, tube-like fullerenes were studied for parameterization. In Section 5, the octahedral, tetrahedral, and D3 fullerenes were parameterized. In section 6, the conclusion and outlook were presented.

## 2. Combination of two vertices

**Theorem 1.** On the fullerene manifold, two adjacent vertices are symmetrically distributed around their parent cone.

**Proof.** Let $A$, $B$ be two adjacent vertices. Cones $A$ and $B$ can be unfolded together as Figure 1(a). First, unfold cone $B$ along a generatrix, and then unfold cone $A$ along line segment $AB$. In the figure, $B_1$, and $B_2$ represent the different positions of vertex $B$ after unfolding.

The angular defect of each vertex is 60°, so the angular defect of their parent vertex is 120°. Let $P$ be the parent vertex of $A\&B$ as Figure (1b), then $PB_1 = PB_2$, $\angle B_1PB_2 = 120°$.

As triangle $AB_1B_2$ is equilateral according to unfolding, point $P$ is its center. Then in Figure 1(c), $\angle APB_1 = \angle APB_2 = 120°$. This means points $A$ and $B$ are symmetrically distributed on the parent cone. Their field angle from the parent vertex is half of the cone angle.

□

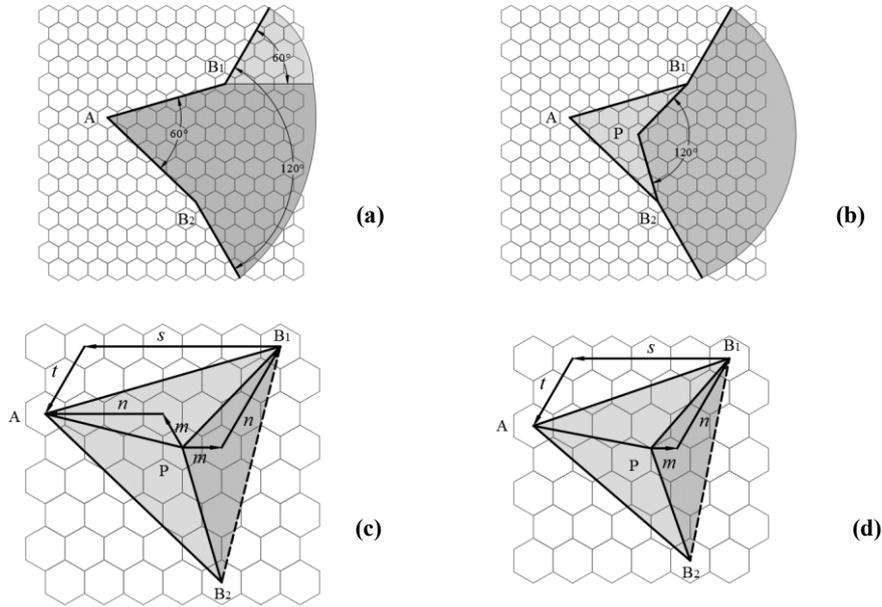

**Fig.1. combination of two vertices.** **(a)** Together unfolding of two adjacent cones $A$ and $B$. **(b)** Position of their parent vertex $P$. **(c)** The relationship of $A$, $B$, and $P$. **(d)** The first kind of unconventional vertex.

The vertices $A$ and $B$ can be obtained by trimming cone $P$: flatten the cone $P$ through line segments $PA$, $PB$, trim off the apex along the connecting line $AB$ and seal the remaining along $AB$, then vertices each with angular defect 60° will appear at points $A$ and $B$.

**Location of the parent vertex**

Because of the symmetry, only a pair of Coxeter coordinates are needed to locate two combined vertices on their parent cone. Let the Coxeter coordinates of $A/B$ to vertex $P$ be $(m, n)$ and the Coxeter coordinates between $A$ and $B$ be $(s, t)$ as Figure1(c), their relationship can be derived with the Eisenstein plane [3] as follows.

The hexagonal lattice is congruent to an Eisenstein plane characterized by a triangular lattice. The Eisenstein plane is a complex plan with unit vectors in the six directions denoted with $1, \omega, \omega^2, \omega^3, \omega^4, \omega^5$, where $\omega = e^{i2\pi/6}$ is the complex root of equation $\omega^6 = 1$. As $\omega^3 = -1$, $\omega^2 = \omega - 1$, each vector in the Eisenstein plane can be express as $(a+b\omega)$, with $(a, b)$ known as Eisenstein integers.

In Figure 2(c), $\overrightarrow{AB_1} = \overrightarrow{AP} + \overrightarrow{PB_1}$, so

$$s + t\omega = (n + m\omega^5) + (m + n\omega)$$
$$= (n + 2m) + (n - m)\omega$$

This means

$$\begin{cases} s = n + 2m \\ t = n - m \end{cases} \quad . \tag{1}$$

Then we can get:

$$\begin{cases} m = (s-t)/3 \\ n = (s-t)/3 + t \end{cases} \quad . \tag{2}$$

**Theorem 2.** If and only if the two Coxeter coordinates between two vertices are congruent modulo 3, their parent vertex is a conventional 4-degree vertex.

**Proof.** According to formula (2), when $s \equiv t \mod (3)$, the Coxeter coordinates $m, n$ are all integers. In this case, the parent vertex corresponds to the center of a square, i.e. a cone vertex with 4 degrees. When $s$ and $t$ are not congruent modulo 3, then $m$ and $n$ are not integers. In this case, the parent vertex corresponds to the position of a carbon atom, just as in figure 1(d). In this article, such vertex is called the first kind of unconventional vertex.

□

## 3. Combination of three vertices

**Theorem 3.** Around the parent cone of three combined vertices, only two vertices are independent.

**Proof.** Let $A$, $B$, and $C$ be three adjacent vertices. Cones $A$, $B$, and $C$ can be unfolded together as Figure 2(a). First, unfold cone $C$ along a generatrix, and then unfold cones $A$ and $B$ along lines segments $CA$, and $CB$ respectively. In the figure, $C_1/C_2/C_2$ represent different positions of vertex C after unfolding, and $C_1$&$C_2$, glued together, represent the location of vertex $C$ on the parent cone.

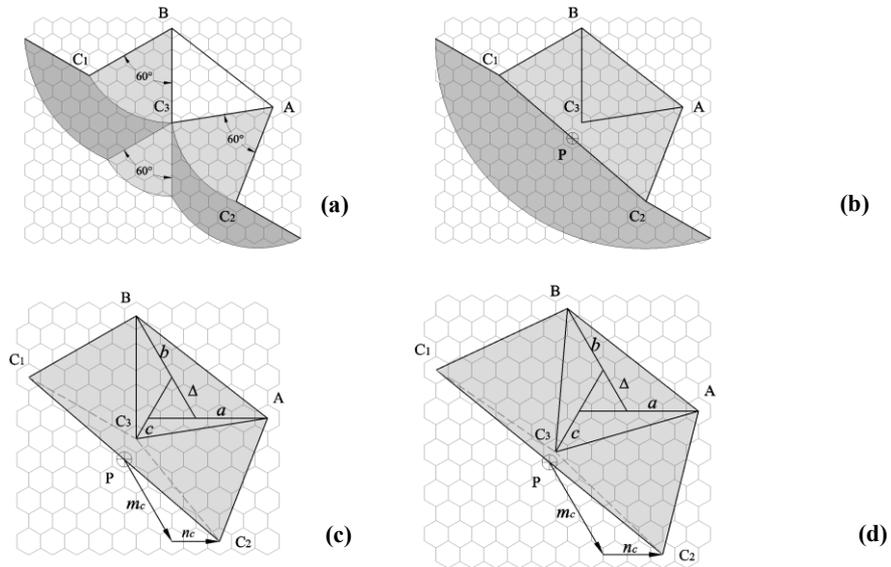

**Fig.2. combination of three vertices.** (a) Together unfolding of three adjacent cones $A$, $B$, and $C$. (b) Position of the parent vertex $P$. (c) The relationship of $A$, $B$, $C$ and $P$. (d) The second kind of unconventional vertex.

Let $P$ be their parent vertex as Figure 2(b). As the total angular defect is 180°, then,

$$PC_1 = PC_2, \angle C_1PC_2 = 180°.$$

This means point P ( marked with ⊕ in the Figure) is just at the midpoint of line segment $C_1C_2$, i.e. points $C_1$ and $C_2$ are symmetrical about point P. As Triangles $BC_1C_3$ and $AC_3C_2$ are equilateral according to unfolding, with the location of C and A we can obtain the position for points $C_1$, $C_2$, and $C_3$, and then for point B. Therefore, only two locations of A, B, and C are independent on their parent cone.

□

According to Theorem 3, only two pairs of Coxeter coordinates are needed to locate three adjacent vertices on their parent cone.

Vertices A, B, and C can be obtained by truncating the parent cone P: take the parent cone as the surface of triangular pyramid P-ABC, truncate off the apex through points A, B, and C, and extend the lattice to section ABC, then vertices each with an angular defect of 60° will appear at points A, B, and C.

**Location of the parent vertex**

For triangle ABC, lattice lines from three vertices form a pinwheel spiral with three arms [4], just as Figure 2(c). Therefore, the triangle ABC has four dimensions, the three spiral arms a, b, c, and the triangle spiral center $\Delta$. Let Coxeter coordinates of C to P be $(m_c, n_c)$, then the Coxeter coordinates between $C_1$ and $C_2$ are $(2m_c, 2n_c)$. With the help of the Eisenstein plane as section 2, the following relationship can be obtained:

$$\begin{cases} m_c = (c+a)/2 + \Delta \\ n_c = (c+b)/2 \end{cases} \quad (3)$$

Similarly, the Coxeter coordinates of B and C can also be obtained with the same method. They can be together written in a formula with matrices：

$$\begin{pmatrix} m_a \\ n_b \\ m_b \\ n_c \\ m_c \\ n_a \end{pmatrix} = \begin{pmatrix} 1/2 & 1/2 & 0 & 1 \\ 1/2 & 1/2 & 0 & 0 \\ 0 & 1/2 & 1/2 & 1 \\ 0 & 1/2 & 1/2 & 0 \\ 1/2 & 0 & 1/2 & 1 \\ 1/2 & 0 & 1/2 & 0 \end{pmatrix} \times \begin{pmatrix} a \\ b \\ c \\ \Delta \end{pmatrix} \quad (4)$$

**Theorem 4.** If and only if the spiral arms of the three adjacent vertices have the same parity, their parent vertex is a conventional 3-degree vertex.

**Proof.** According to formula (4), if a, b, and c have the same parity, the Coxeter coordinates of points A, B, C to point P are all integers. In this case, the parent vertex corresponds to the center of a triangle, i.e. a cone vertex with three degrees.

If a, b, and c don't have the same parity, the Coxeter coordinates of points A, B, C to point P are not all integers. In this case, the parent vertex corresponds to the center of a carbon bond on the graphene sheet, just as Figure 2(d). In this paper, such a parent vertex is called the second kind of unconventional vertex.

□

## 4. Parameterization of tube-like fullerenes

Most fullerenes have a tube-like structure as carbon nanotubes. The twelve vertices of a tube-like fullerene cluster into two sextuples, each forming a cap at the end of the tube. The cap has five possible shapes according to the number of vertices adjacent to the tube: hexagon, pentagonal pyramid, trimmed square pyramid, truncated triangular pyramid, and a double trimmed digonal pyramid.

**Parameters for a fullerene capped with two pentagonal pyramids**

If only five vertices are adjacent to the tube, the cap is a pentagonal pyramid. Five pairs of Coxeter coordinates of the base edges can be selected as its independent parameters. Just as Figure 3(a), as the triangle $AB_1B_2$ is equilateral according to unfolding, the position of vertex $A$ is determined, and then the detail of the cap can be obtained.

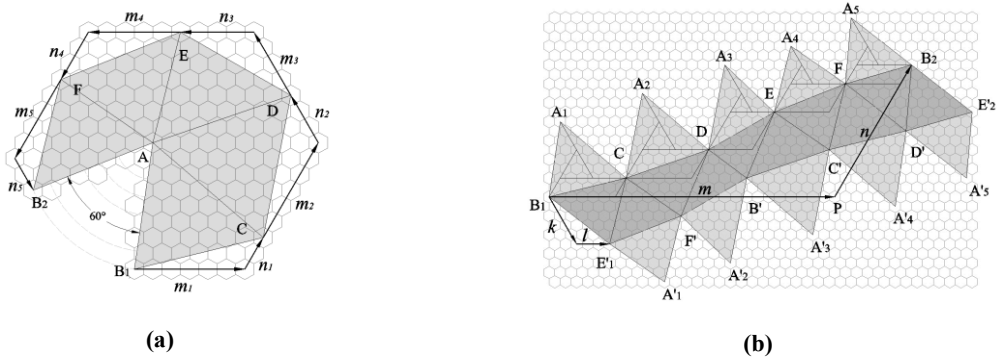

**Fig.3. Parameters for a fullerene capped with two pentagonal pyramids.** (a) Parameters for a cap as a pentagonal pyramid. (b) A fullerene capped with two pentagonal pyramids. Because two caps are connected with a tube, there are two constrain equations for the parameters of these two caps. Two other parameters are needed for their relative position.

A tube-like fullerene capped with two pentagonal pyramids can be unfolded as Figure 3(b). Let the feature parameters of the tube be $(m, n)$. To match the tube, the parameters of a pentagonal pyramid must satisfy the following two equations:

$$m = \sum_{i=1}^{5} m_i \quad ; \quad n = \sum_{i=1}^{5} n_i \; .$$

Let $(m_i', n_i')$ be the parameters of the opposite pentagonal pyramid. then we can get:

$$\begin{cases} \sum_{i=1}^{5} m_i = \sum_{j=1}^{5} m_j' \\ \sum_{i=1}^{5} n_i = \sum_{j=1}^{5} n_j' \end{cases} \quad . \tag{5}$$

According to formula (5), only 18 of 20 parameters of the two caps are independent. Additionally, another two parameters are needed to determine their relative position, just as $(k, l)$ in Figure 3(b). Together, they make up 20 independent parameters for the fullerene.

**Parameters for other caps**

If only four vertices are adjacent to the tube, the cap is a trimmed square pyramid. Just as in figure 4(a), four vertices $A$, $B$, $C$, and $D$ are adjacent to the tube. They and the parent vertex $P$ of $E\&F$ form a square pyramid $P$-$ABCD$. Four pairs of Coxeter coordinates of the base edges can be selected as the independent parameters for the pyramid. According to section 2, another pair of Coxeter coordinates are needed to locate vertices $E$ and $F$. To avoid possible non-integer values, the Coxeter coordinates for $E$ and $F$ can refer to a base point of the pyramid instead of their parent vertex, just as $(s, t)$ in the Figure. This selection method is also suitable for other caps.

If only three vertices are adjacent to the tube, the cap is a truncated triangular pyramid. Just as Figure 4(b), three vertices *A*, *B*, *C* are adjacent to the tube. They and the parent vertex *P* of *D*, *E*, and *F* form a triangular pyramid *P-ABC*. Three pairs of Coxeter coordinates of the base edges can be selected as the independent parameters for the triangular pyramid. According to section 3, another two pairs of Coxeter coordinates are needed for vertices *D*, *E*, and *F*.

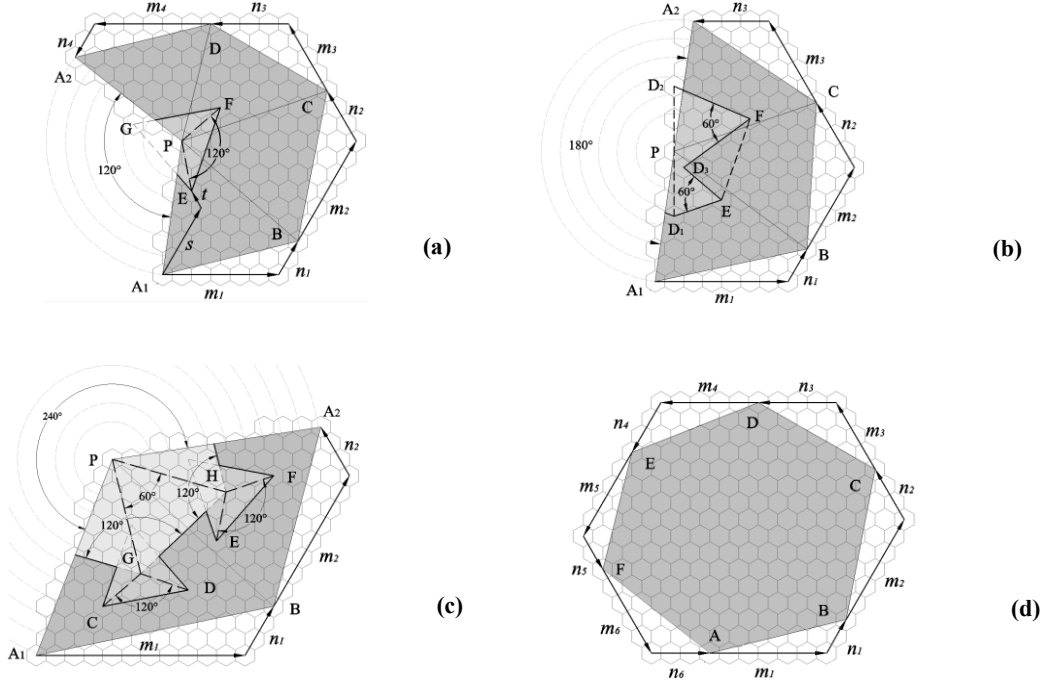

**Fig.4. Other kinds of caps. (a)** Trimmed square pyramid. **(b)** Truncated triangular pyramid. **(c)** Double trimmed dihedral pyramid. **(d)** Hexagonal cap.

If only two vertices are adjacent to the tube, the cap is a double trimmed digonal pyramid. Just as Figure 4(c), vertices *A* and *B* are adjacent to the tube. *G* is the parent vertex of *C&D*. *H* is the parent vertex of *E&F*. Vertices *G* and *H* can also be combined to get a parent vertex *P*, which is the apex of the digonal pyramid *P-AB* that has only two lateral faces.

Two pairs of Coxeter coordinates for the base edges can be selected as independent parameters for the digonal pyramid. With the same angular defeat, vertices *G* and *H* are symmetrically distributed on their parent cone *P*. This means that ∠*GPH* = 60°, and *GHP* is equilateral triangle. With the symmetries of *G&H*, *C&D*, and *E&F*, only three pair of Coxeter coordinates are needed for the locations of vertices *C*, *D*, *E*, and *F* on the pyramid.

If six vertices of a cap are all adjacent to the tube, they form a hexagon as Figure 4(d). As a close loop, only five of its six edges are independent. Their Coxeter coordinates can be selected as independent parameters for the cap. With the help of the Eisenstein plane, we can get the Coxeter coordinates of the sixth edge from these parameters:

$$\begin{cases} m_6 = m_2 + m_3 - m_5 + n_1 + n_2 - n_4 - n_5 \\ n_6 = -m_1 - m_2 + m_4 + m_5 - n_1 + n_3 + n_4 \end{cases}$$

**Parameterization of any tube-like fullerene**

A hexagonal cap can be treated as a special hexagonal pyramid with zero angular defect at the apex. Therefore, each tube-like fullerene has a parent structure as a tube capped with two polygonal pyramids. To match the common tube in the middle, the parameters for the two polyhedral pyramids have two constrain equations similar as formula (5), with ranges of *i* and *j* according to the base number of the pyramids. Each cap has 10 independent parameters. Only 18 of these 20 parameters are independent. Together with two parameters for the relative position of the two caps, as $(k, l)$ in figure 3(b), they make up 20 independent parameters for a tube-like fullerene.

## 5. Parameterization of octahedral fullerenes, tetrahedral fullerenes, and $D_3$-fullerenes

Except for clustering into two sextuples as a tube-like fullerene, the 12 vertices of a fullerene may also evenly cluster into six pairs, four triples, or three quads. Such fullerenes are called octahedral fullerenes, tetrahedral fullerenes, and $D_3$-fullerenes respectively.

**Parameters for octahedral fullerenes**

An octahedral fullerene is a trimmed (4,6)-fullerene if the parent vertex of each pair is a conventional 4-degree vertex. As each pair of vertices have two independent parameters according to section 2, a trimmed (4,6)-fullerene has 12 such parameters in total. Together with the 8 parameters for (4,6)-fullerene in literature [4], they make up 20 independent parameters for the trimmed (4,6)-fullerene.

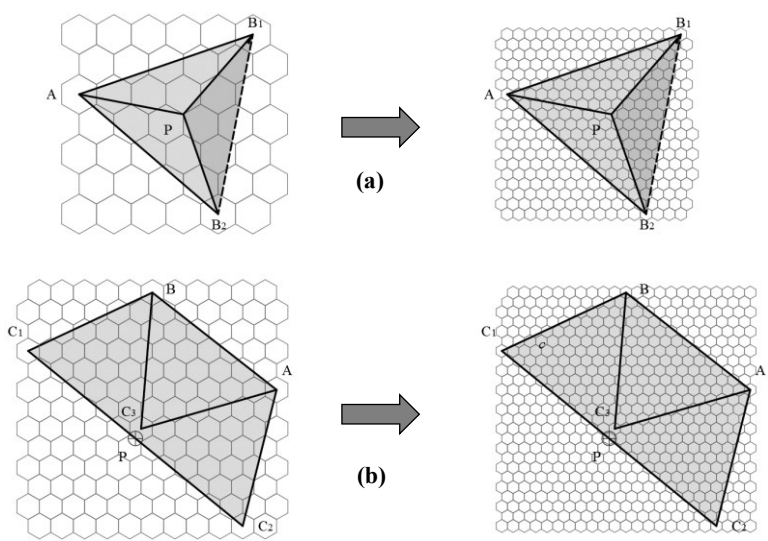

**Fig.5. Transformation of unconventional vertices.** (a) If the lattice is refined with one-third-sized hexagons, the first kind of unconventional vertex will be transformed to a conventional 4-degree vertex. (b) If the lattice is refined with half-sized hexagons, the second kind of unconventional vertex will be transformed to a conventional 3-degree vertex.

The parameters for the trimmed (4,6)-fullerenes are suitable for all the octahedral fullerenes. Just as in Figure 5(a), with lattice properly refined, the first kind of unconventional vertex can be turned into a conventional 4-degree vertex, and any other octahedral fullerene can be transformed to a trimmed (4,6)-fullerene. Therefore, they can be parameterized as a trimmed (4,6)-fullerene. The difference is that the values of these parameters may not be integers.

**Parameters for tetrahedral fullerenes**

A tetrahedral fullerene is a truncated (3,6)-fullerene if the parent vertex for each triple is a conventional 3-degree vertex. As each triple of vertices has four independent parameters according to section 3, a truncated (3,6)-fullerene has 16 such parameters in total. Together with the 4 independent parameters for (3,6)-fullerene in literature [4], they make up 20 independent parameters for the truncated (3,6)-fullerene.

Similarly, the second kind of unconventional vertices can be turn into conventional 3-degree vertices with lattice refined as Figure 5(b). Therefore, other tetrahedral fullerenes can also be parameterized as a truncated (3,6)-fullerene, with parameter values may not be integers.

**Parameters for $D_3$-fullerenes**

As the parent structure is a triangular dihedron, a fullerene whose vertices cluster into three quads is called a $D_3$-fullerene. The triangular dihedron has only two independent parameters [2]. As each quad has six independent parameters as mentioned in section 4, a $D_3$-fullerene has 18 such parameters in total. Together with the two parameters for the triangular dihedron, they make up 20 independent parameters for the $D_3$-fullerene.

## 6. Conclusion and outlook

Parameterization with vertex combinations is a practical method for many types of fullerenes. With both the information of parent structure as framework and the symmetries of local clusters, such parameters may be different for fullerenes from type to type, but they are able to characterize the specific shape for each type.

A pair of Coxeter coordinates correspond to two possible trimming results for a four-degree vertex. Two pairs of Coxeter coordinates correspond to up to six possible truncations for a three-degree vertex. Therefore, it is essential to specify the Coxeter coordinates for parameters, or to eliminate ambiguity with further classification.

**References**

[1] W.P. Thurston, Shapes of polyhedra and triangulations of the sphere, Geom. Top. Monographs, 1 (1998) 511–549.

[2] P.W. Fowler, D.E. Manolopoulos, An atlas of fullerenes, 2nd ed. Mineola, NY: Dover Publications Inc.2006.

[3] P. Schwerdtfeger, L.N. Wirz, J. Avery, The topology of fullerenes, Wiley Interdiscip. Rev. Comput. Mol. Sci. 5 (2015) 96–145.

[4] S. Li, On the intrinsic symmetries and parameterization of fullerenes, Comput. Theor. Chem. 1200 (2021).

[5] E.J. Hartung, The Clar structure of fullerenes, Syracuse University, 2012, pp. 112.